\documentclass[11pt]{article}

\usepackage{amssymb,latexsym}
\usepackage{epsfig}
\usepackage{eufrak}
\usepackage{amsmath}
\usepackage{mathrsfs}
\usepackage{color}

\newtheorem{theorem}{Theorem}
\newtheorem{lemma}{Lemma}
\newtheorem{corollary}{Corollary}

\newcommand{\be}{\begin{equation}}
\newcommand{\ee}{\end{equation}}
\newcommand{\bee}{\begin{eqnarray*}}
\newcommand{\eee}{\end{eqnarray*}}
\newcommand{\bel}{\begin{eqnarray}}
\newcommand{\eel}{\end{eqnarray}}
\newcommand{\bec}{\begin{cases}}
\newcommand{\eec}{\end{cases}}
\newcommand{\bem}{\begin{bmatrix}}
\newcommand{\eem}{\end{bmatrix}}

\newcommand{\la}{\label}
\newcommand{\li}{\left}
\newcommand{\ri}{\right}

\newcommand{\lc}{\lceil}
\newcommand{\rc}{\rceil}
\newcommand{\lf}{\lfloor}
\newcommand{\rf}{\rfloor}

\newcommand{\vep}{\varepsilon}

\newcommand{\de}{\delta}

\newcommand{\ze}{\zeta}
\newcommand{\al}{\alpha}

\newcommand{\ro}{\rho}

\newcommand{\f}{\frac}

\newcommand{\cd}{\cdots}

\newcommand{\qu}{\quad}
\newcommand{\qqu}{\qquad}
\newcommand{\fa}{\forall}

\newcommand{\mscr}{\mathscr}
\newcommand{\mcal}{\mathcal}
\newcommand{\mbf}{\mathbf}

\newcommand{\wh}{\widehat}
\newcommand{\wt}{\widetilde}
\newcommand{\mrm}{\mathrm}
\newcommand{\bs}{\boldsymbol}

\newcommand{\sh}{\slash}

\newcommand{\tx}{\text}

\newcommand{\bed}{\begin{description}}
\newcommand{\eed}{\end{description}}
\newcommand{\bei}{\begin{itemize}}
\newcommand{\eei}{\end{itemize}}
\newcommand{\ben}{\begin{enumerate}}
\newcommand{\een}{\end{enumerate}}
\newcommand{\bib}{\bibitem}
\newcommand{\beL}{\begin{lemma}}
\newcommand{\eeL}{\end{lemma}}
\newcommand{\beT}{\begin{theorem}}
\newcommand{\eeT}{\end{theorem}}
\newcommand{\sect}{\section}

\newcommand{\bpf}{\begin{pf}}
\newcommand{\epf}{\end{pf}}
\newcommand{\bsk}{\bigskip}
\newcommand{\bi}{\binom}

\newcommand{\pfbox}{\hfill\mbox{$\Box$}}

\newenvironment{pf}{\paragraph*{Proof{\rm.}}}{\pfbox\bigskip}

\begin{document}

\title{{\bf On Estimation of Finite Population Proportion}
\thanks{The author is currently with Department of Electrical Engineering,
Louisiana State University at Baton Rouge, LA 70803, USA, and
Department of Electrical Engineering, Southern University and A\&M College, Baton
Rouge, LA 70813, USA; Email: chenxinjia@gmail.com}}

\author{Xinjia Chen}


\maketitle

\begin{abstract}

In this paper, we study the classical problem of estimating the
proportion of a finite population.   First, we consider a fixed
sample size method and derive an explicit sample size formula which
ensures a mixed criterion of absolute and relative errors.   Second,
we consider an inverse sampling scheme such that the sampling is
continue until the number of units having a certain attribute
reaches a threshold value or the whole population is examined.  We
have established a simple method to determine the threshold so that
a prescribed relative precision is guaranteed.  Finally, we develop
a multistage  sampling scheme for constructing fixed-width
confidence interval for the proportion of a finite population.
Powerful computational techniques are introduced  to make it
possible that the fixed-width confidence interval ensures prescribed
level of coverage probability.

\end{abstract}

\section{Fixed Sample Size Method}

The estimation of the proportion of a finite population is a basic
and very important problem in probability and statistics \cite{Desu,
Thompson}.  Such problem finds applications spanning many areas of
sciences and engineering.  The problem is formulated as follows.

Consider a finite population of $N$ units, among which there are $M$
units having a certain attribute. The objective is to estimate the
proportion $p = \f{M}{N}$ based on sampling without replacement.

One popular method of sampling is to draw $n$ units without
replacement from the population and count the number, $\mathbf{k}$,
of units having the attribute.  Then, the estimate of the proportion
is taken as $\wh{\bs{p}} = \f{\mathbf{k}}{n}$. In this process, the
sample size $n$ is fixed.

Clearly, the random variable $\mathbf{k}$ possesses a hypergeometric
distribution. The reliability of the estimator $\wh{\bs{p}} =
\f{\mathbf{k}}{n}$ depends on $n$.  For error control purpose, we
are interested in a crucial question as follows:

\bsk

{\it For prescribed  margin of absolute error $\vep_a \in (0, 1)$,
margin of relative error $\vep_r \in (0, 1)$, and confidence
parameter $\de \in (0, 1)$, how large the sample size $n$ should be
to guarantee \be
 \la{cov}
 \Pr \li \{ |\wh{\bs{p}} - p| < \vep_a  \; \; \mrm{or} \;\; \li | \f{
\wh{\bs{p}} - p } {p} \ri | < \vep_r  \ri \} > 1 - \de ? \ee }

In this regard, we have

\beT Let $\vep_a \in (0, 1)$ and $\vep_r \in (0, 1)$ be real numbers
such that $\f{ \vep_a}{ \vep_r} + \vep_a \leq \f{1}{2}$.   Then,
(\ref{cov}) is guaranteed provided that \be \la{con}
 n > \f{ \vep_r \ln \f{2}{\de} } { \li ( \vep_a + \vep_a \vep_r \ri ) \ln (1 + \vep_r) + \li ( \vep_r - \vep_a -
\vep_a \vep_r \ri ) \ln \li ( 1 -  \f{ \vep_a \vep_r } { \vep_r - \vep_a  }
 \ri )  }.
\ee \eeT

\bsk

The proof of Theorem 1 is given in Appendix A. It should be noted
that conventional methods for determining sample sizes are based on
normal approximation, see \cite{Desu} and the references therein. In
contrast, Theorem 1 offers a rigorous method for determining sample
sizes.  To reduce conservativeness, a numerical approach has been
developed by Chen \cite{Chen2} which permits exact computation of
the minimum sample size.

\section{Inverse Sampling of Finite Population}

To estimate the proportion $p$, a frequently-used sampling method is
the {\it inverse sampling} scheme described as follows:

Continuing sampling from the population (without replacement) until
$r$ units found to carry the attribute or the number of sample size
$\bs{n}$ reaches the population size $N$.  The estimator of the
proportion $p$ is taken as the ratio $\wt{\bs{p}} =
\f{\bs{k}}{\bs{n}}$, where $\bs{k}$ is the number of units having
the attribute among the $\bs{n}$ units.

Clearly, the reliability of the estimator $\wt{\bs{p}}$ depends on
the threshold value $r$.  Hence, we are interested in a crucial
question as follows:

\bsk

{\it For prescribed margin of relative error $\vep \in (0, 1)$ and
confidence parameter $\de \in (0, 1)$, how large the threshold $r$
should be to guarantee
 \[
 \Pr \li \{  \li |
\wt{\bs{p}} - p  \ri | < \vep p \ri \} > 1 - \de ? \] }

For this purpose, we have

 \beT \la{main_Binomial} For any $\vep \in (0,1)$,
\[
\Pr \li \{  \li |  \wt{\bs{p}} - p  \ri | \geq \vep p \ri \} \leq
\mscr{Q} (\vep, r) \]
 where {\small
 \[
 \mscr{Q} (\vep, r) =
(1 + \vep)^{- r} \exp \li ( \f{ \vep r }{ 1 + \vep } \ri ) + (1 -
\vep)^{- r} \exp \li ( - \f{ \vep r }{ 1 - \vep }  \ri ), \] } which
is monotonically decreasing with respect to $r$. Moreover, for any
$\de \in (0,1)$, there exists a unique number $r^*$ such that
$\mscr{Q} (\vep, r^*) = \de$ and {\small
\[
\max \li \{ \f{ (1 + \vep) \ln \f{1}{\de} } { (1 + \vep) \ln (1 +
\vep) - \vep }, \; \f{ (1 - \vep) \ln \f{2}{\de} } { (1 - \vep) \ln
(1 - \vep) + \vep } \ri \} < r^* < \f{ (1 + \vep) \ln \f{2}{\de} } {
(1 + \vep) \ln (1 + \vep) - \vep }. \] } \eeT

\bsk

The proof of Theorem 2 is given in Appendix B.  As an immediate
consequence of Theorem 2, we have

\begin{corollary}
\la{exp}
 Let $\vep, \; \de \in (0,1)$.  Then, $ \Pr \li \{  \li |  \wt{\bs{p}} - p \ri | <
\vep p \ri \} > 1 - \de$ provided that \be \la{formu}
 \boxed{r > \f{ (1 + \vep) \ln \f{2}{\de}  } { (1 + \vep) \ln (1 + \vep) - \vep
} } \ee
\end{corollary}

\sect{Multistage Fixed-width Confidence Intervals}

So far we have only considered point estimation for the proportion
$p$.  Interval estimation is also an important method for estimating
$p$.  Motivated by the fact that a confidence interval must be
sufficiently narrow to be useful, we shall develop a multistage
sampling scheme for constructing a fixed-width confidence interval
for the proportion, $p$, of the finite population discussed in
previous sections.

Note that the procedure of sampling without replacement can be
precisely described as follows:

 Each time a single unit is drawn without replacement from the remaining population so
that every unit of the remaining population has equal chance of
being selected.

Such a sampling process can be exactly characterized by random
variables $X_1, \cd, X_N$ defined in a probability space $(\Omega,
\mscr{F}, \Pr)$ such that $X_i$ denotes the characteristics of the
$i$-th sample in the sense that $X_i = 1$ if the $i$-th sample has
the attribute and $X_i = 0$ otherwise.  By the nature of the
sampling procedure, it can be shown that
\[ \Pr \{ X_i = x_i, \; i = 1, \cd, n \} =
\bi{M}{\sum_{i = 1}^n x_i} \bi{N - M}{n - \sum_{i = 1}^n x_i}
 \li \slash \li [ \bi{n}{\sum_{i = 1}^n x_i} \bi{N}{n} \ri. \ri ]
\]
for any $n \in \{1, \cd, N\}$ and any $x_i \in \{0, 1\}, \; i = 1,
\cd, n$.  Based on random variables $X_1, \cd, X_N$, we can define a
multistage sampling scheme of the following basic structure.  The
sampling process is divided into $s$ stages with sample sizes $n_1 <
n_2 < \cd < n_s$. The continuation or termination of sampling is
determined by decision variables.  For each stage with index $\ell$,
a decision variable $\bs{D}_\ell = \mscr{D}_\ell (X_1, \cd,
X_{n_\ell})$ is defined based on random variables $X_1, \cd,
X_{n_\ell}$. The decision variable $\bs{D}_\ell$ assumes only two
possible values $0, \; 1$ with the notion that the sampling is
continued until $\bs{D}_\ell = 1$ for some $\ell \in \{1, \cd, s\}$.
Since the sampling must be terminated at or before the $s$-th stage,
it is required that $\bs{D}_s = 1$. For simplicity of notations, we
also define $\bs{D}_\ell = 0$ for $\ell = 0$.

Our goal is to construct a fixed-width confidence interval $(\bs{L},
\bs{U})$ such that $\bs{U} - \bs{L} \leq 2 \vep$ and that $\Pr \{
\bs{L} < p < \bs{U} \mid p  \}
> 1 - \de$ for any $p \in \{ \f{i}{N}: 0 \leq i \leq N \}$ with prescribed $\vep \in (0, \f{1}{2})$ and $\de \in (0,
1)$.  Toward this goal, we need to define some multivariate
functions as follows.

For $\al \in (0, 1)$ and integers $0 \leq k \leq n \leq N$, let
$\mcal{L} (N, n, k, \al)$ be the smallest integer $M_l$ such that
{\small $\sum_{i = k }^{n} \bi{M_l}{i} \bi{N - M_l}{n - i} \sh
\bi{N}{n} > \f{\al}{2}$}.  Let $\mcal{U} (N, n, k, \al)$ be the
largest integer $M_u$ such that {\small $\sum_{i = 0}^{k}
\bi{M_u}{i} \bi{N - M_u}{n - i} \sh \bi{N}{n}
> \f{\al}{2}$}.  Let $n_{\mrm{max}} (N, \al)$ be the smallest
number $n$  such that $\mcal{U} (N, n, k, \al) - \mcal{L} (N, n, k,
\al) \leq 2 \vep N$ for $0 \leq k \leq n$. Let $n_{\mrm{min}} (N,
\al)$ be the largest number $n$  such that $\mcal{U} (N, n, k, \al)
- \mcal{L} (N, n, k, \al) > 2 \vep N$ for $0 \leq k \leq n$.

\beT \la{coverage_abs} Let $\ze > 0$ and $\ro > 0$. Let $n_1 < n_2 <
\cd < n_s$ be the ascending arrangement of all distinct elements of
{\small $ \li \{ \li \lc \li [ \f{n_{\mrm{max}} (N, \ze \de)}{
n_{\mrm{min}} (N, \ze \de)}  \ri ]^{ \f{i}{\tau} } n_{\mrm{min}} (N,
\ze \de)  \ri \rc : i = 0, 1, \cd, \tau \ri \}$} with {\small  $\tau
= \li \lc \f{1}{ \ln (1 + \ro)}  \ln \f{n_{\mrm{max}} (N, \ze \de)}{
n_{\mrm{min}} (N, \ze \de)}  \ri \rc$. }  For $\ell = 1, \cd, s$,
define $K_\ell = \sum_{i = 1}^{n_\ell} X_i$ and $\bs{D}_\ell$ such
that $\bs{D}_\ell = 1$ if $\mcal{U} (N, n_\ell, K_\ell, \ze \de) -
\mcal{L} (N, n_\ell, K_\ell, \ze \de) \leq  2 \vep N$; and
$\bs{D}_\ell = 0$ otherwise. Suppose the stopping rule is that
sampling is continued until $\bs{D}_\ell = 1$ for some $\ell \in
\{1, \cd, s\}$.  Define {\small $\bs{L} = \f{1}{N} \times \mcal{L}
\li (N, \mathbf{n}, \sum_{i=1}^{\mathbf{n}} X_i, \ze \de \ri )$} and
{\small $\bs{U} = \f{1}{N} \times \mcal{U} \li (N, \mathbf{n},
\sum_{i=1}^{\mathbf{n}} X_i, \ze \de \ri )$}, where $\mathbf{n}$ is
the sample size when the sampling is terminated. Then, a sufficient
condition to guarantee $\Pr \li \{ \bs{L} < p < \bs{U} \mid p \ri \}
> 1 - \de$ for any $p \in \{ \f{i}{N}: 0 \leq i \leq N \}$ is that {\small  \bel &  &
\sum_{\ell = 1}^s \li [ \Pr \{ \mcal{L} (N, n_\ell, K_\ell, \ze \de)
\geq
M, \; \bs{D}_{\ell - 1} = 0, \; \bs{D}_\ell = 1 \mid M \} \ri. \nonumber\\
&  &  + \li .  \Pr \{ \mcal{U} (N, n_\ell, K_\ell, \ze \de) \leq M,
\; \bs{D}_{\ell - 1} = 0, \; \bs{D}_\ell = 1 \mid M \}  \ri ] < \de
\la{2D2}
 \eel}
 for all $M \in \{0, 1, \cd, N \}$,
where (\ref{2D2}) is satisfied if  $\ze > 0$ is sufficiently small.

\eeT

\bsk

It should be noted that Theorem 3 has employed the
double-decision-variable method recently proposed by Chen in
\cite{Chen_EST}.  To further reduce computational complexity, the
techniques of bisection confidence tuning and domain truncation
developed in \cite{Chen_EST, ChenT} can be very useful.

\appendix

\sect{Proof of Theorem 1}

To prove the theorem, we shall introduce function
\[
g(\vep, p) =  (p + \vep) \ln \f{ p}{p + \vep} + (1 - p - \vep) \ln
\f{1 - p}{1 - p - \vep}
\]
where $0 < \vep < 1 - p$.   We need some preliminary results.

The following lemma is due to Hoeffding \cite{Hoeffding}.

\beL \la{lem1}
\[
\Pr \{ \wh{\bs{p}} \geq p + \vep \} \leq \exp ( n \; g(\vep, p) )
\qu \tx{for} \qu 0 < \vep < 1 - p < 1,
\]
\[
\Pr \{ \wh{\bs{p}} \leq p - \vep \} \leq \exp ( n \; g(-\vep, p) )
\qu \tx{for} \qu 0 < \vep < p < 1.
\]
\eeL

The following Lemmas 2--4 have been established in \cite{Chen1}.

\beL  \la{lem2}

Let $0 < \vep < \f{1}{2}$. Then, $g(\vep, p)$ is monotonically
increasing with respective to $p \in (0,  \f{1}{2} - \vep)$ and
monotonically decreasing with respective to $p \in (\f{1}{2},  1 -
\vep)$.  Similarly, $g(- \vep, p)$ is monotonically increasing with
respective to $p \in (\vep, \f{1}{2})$ and monotonically decreasing
with respective to $p \in (\f{1}{2} + \vep,  1)$.\eeL

\beL \la{lem3}

Let $0 < \vep < \f{1}{2}$.  Then,
\[
g(\vep, p) > g(-\vep, p) \qqu \fa p \in \li (\vep, \f{1}{2} \ri ],
\]
\[
g(\vep, p) < g(-\vep, p) \qqu \fa p \in \li (\f{1}{2}, 1 - \vep \ri
).
\]
\eeL

\beL \la{lem4}

Let $0 < \vep < 1$.  Then, $g \li ( \vep p, p \ri )$ is
monotonically decreasing with respect to {\small $p \in \li (0,
\f{1}{1 + \vep} \ri )$}. Similarly, $g \li (- \vep p, p \ri )$ is
monotonically decreasing with respect to $p \in (0, 1)$. \eeL

\beL \la{lem5}  Suppose $0 < \vep_r < 1$ and $0 < \f{\vep_a}{\vep_r} + \vep_a \leq \f{1}{2}$.  Then, \be \la{ineqm}
 \Pr \{ \wh{\bs{p}} \leq p - \vep_a \} \leq \exp \li (n \; g \li ( -
\vep_a, \f{\vep_a}{\vep_r}\ri ) \ri ) \ee for $0 < p \leq
\f{\vep_a}{\vep_r}$.  \eeL

\bpf We shall show (\ref{ineqm}) by investigating three cases as
follows.  In the case of $p < \vep_a$, it is clear that
\[
\Pr \{ \wh{\bs{p}} \leq p - \vep_a \} = 0  < \exp \li (n \; g \li (
- \vep_a, \f{\vep_a}{\vep_r}\ri ) \ri ).
\]

In the case of $p = \vep_a$, we have \bee \Pr \{ \wh{\bs{p}} \leq p - \vep_a \} & = &
\Pr \{ \wh{\bs{p}} = 0  \} = \Pr \{ \mbf{k} = 0 \}\\
&  = & \f{ \bi{N - M}{n} } { \bi{N}{n} } \leq \li ( \f{N - M}{N} \ri )^n \\
& = & (1 - p)^n =  (1 - \vep_a)^n\\
& = & \lim_{p \to \vep_a} \exp (n \; g(- \vep_a, p))\\
& < &  \exp \li (n \; g \li ( - \vep_a, \f{\vep_a}{\vep_r}\ri ) \ri ), \eee where the last inequality follows from Lemma \ref{lem2} and the fact
that $\vep_a < \f{\vep_a}{\vep_r} \leq \f{1}{2} - \vep_a$.

In the case of $\vep_a < p \leq \f{\vep_a}{\vep_r}$,  we have
\[
\Pr \{ \wh{\bs{p}} \leq p - \vep_a \} \leq \exp (n \; g(- \vep_a,
p)) < \exp \li (n \; g \li ( - \vep_a, \f{\vep_a}{\vep_r}\ri ) \ri
),
\]
where the first inequality follows from Lemma \ref{lem1} and the second inequality follows from Lemma \ref{lem2} and the fact that $\vep_a <
\f{\vep_a}{\vep_r} \leq \f{1}{2} - \vep_a$.  So, (\ref{ineqm}) is established. \epf

\beL \la{lem6}  Suppose $0 < \vep_r < 1$ and $0 < \f{\vep_a}{\vep_r} + \vep_a \leq \f{1}{2}$.   Then, \be \la{ineqmm}
 \Pr \{ \wh{\bs{p}} \geq (1 + \vep_r) p \} \leq \exp \li (n \; g \li (
\vep_a, \f{\vep_a}{\vep_r}\ri ) \ri ) \ee for $\f{\vep_a}{\vep_r} <
p < 1$.  \eeL

\bpf We shall show (\ref{ineqmm}) by investigating three cases as
follows.  In the case of $p > \f{1}{1 + \vep_r}$, it is clear that
\[
\Pr \{ \wh{\bs{p}} \geq (1 + \vep_r) p \} = 0  < \exp \li (n \; g
\li ( \vep_a, \f{\vep_a}{\vep_r}\ri ) \ri ).
\]

In the case of $p = \f{1}{1 + \vep_r}$, we have \bee \Pr \{ \wh{\bs{p}} \geq (1 + \vep_r) p \} & = &
\Pr \{ \wh{\bs{p}} = 1  \} = \Pr \{ \mbf{k} = n \}\\
&  = & \f{ \bi{M}{n} } { \bi{N}{n} } \leq \li ( \f{M}{N} \ri )^n \\
& = & p^n =  \li (\f{1}{1 + \vep_r} \ri )^n\\
& = & \lim_{p \to \f{1}{1 + \vep_r}} \exp (n \; g( \vep_r p, p))\\
& < &  \exp \li (n \; g \li ( \vep_a, \f{\vep_a}{\vep_r}\ri ) \ri ), \eee where the last inequality follows from Lemma \ref{lem4} and the fact
that $ \f{\vep_a}{\vep_r} \leq \f{1}{2} \f{1}{1 + \vep_r} < \f{1}{1 + \vep_r}$ as a result of $0 < \f{\vep_a}{\vep_r} + \vep_a \leq \f{1}{2}$.

In the case of $\f{\vep_a}{\vep_r} < p < \f{1}{1 + \vep_r}$,  we
have
\[
\Pr \{ \wh{\bs{p}} \leq (1 + \vep_r) p \} \leq \exp (n \; g(\vep_r
p, p)) < \exp \li (n \; g \li ( \vep_a, \f{\vep_a}{\vep_r}\ri ) \ri
),
\]
where the first inequality follows from Lemma \ref{lem1} and the second inequality follows from Lemma \ref{lem4}.  So, (\ref{ineqmm}) is
established. \epf

\bsk

We are now in a position to prove the theorem. We shall assume (\ref{con}) is satisfied and show that (\ref{cov}) is true.  It suffices to show
that
\[
\Pr \{ |\wh{\bs{p}} - p | \geq \vep_a, \; |\wh{\bs{p}} - p | \geq
\vep_r p \} < \de.
\]

For $0 < p \leq \f{\vep_a}{\vep_r}$, we have \bel \Pr \{
|\wh{\bs{p}} - p | \geq \vep_a, \; |\wh{\bs{p}} - p | \geq \vep_r p
\} & =
& \Pr \{ |\wh{\bs{p}} - p | \geq \vep_a  \} \nonumber\\
& = & \Pr \{ \wh{\bs{p}}  \geq p + \vep_a \} + \Pr \{ \wh{\bs{p}}
\leq p - \vep_a \} \la{ineq6}. \eel

Noting that $0 < p + \vep_a \leq \f{\vep_a}{\vep_r} + \vep_a \leq
\f{1}{2}$, we have
\[
\Pr \{ \wh{\bs{p}}  \geq p + \vep_a \} \leq \exp (n \; g(\vep_a, p))
\leq \exp \li (n \; g \li ( \vep_a, \f{\vep_a}{\vep_r} \ri ) \ri ),
\]
where the first inequality follows from Lemma \ref{lem1} and the second inequality follows from Lemma \ref{lem2}.  It can be checked that
(\ref{con}) is equivalent to
\[
\exp \li (n \; g \li ( \vep_a, \f{\vep_a}{\vep_r} \ri ) \ri ) < \f{\de}{2}.
\]
Therefore,
\[
\Pr \{ \wh{\bs{p}}  \geq p + \vep_a \}  < \f{\de}{2}
\]
for $0 < p \leq \f{\vep_a}{\vep_r}$. \bsk

On the other hand, since $\vep_a < \f{\vep_a}{\vep_r} < \f{1}{2}$, by Lemma \ref{lem5} and Lemma \ref{lem3}, we have
\[
 \Pr \{ \wh{\bs{p}} \leq p - \vep_a \} \leq \exp \li (n \; g \li ( -
\vep_a, \f{\vep_a}{\vep_r}\ri ) \ri ) \leq \exp \li (n \; g \li ( \vep_a, \f{\vep_a}{\vep_r}\ri ) \ri ) < \f{\de}{2}
\]
for $0 < p \leq \f{\vep_a}{\vep_r}$.  Hence, by (\ref{ineq6}),
\[
\Pr \{ |\wh{\bs{p}} - p | \geq \vep_a, \; |\wh{\bs{p}} - p | \geq
\vep_r p \} < \f{\de}{2} + \f{\de}{2} = \de.
\]
This proves (\ref{cov}) for $0 < p \leq \f{\vep_a}{\vep_r}$.

\bsk

For $\f{\vep_a}{\vep_r} < p < 1$, we have \bee \Pr \{ |\wh{\bs{p}} -
p | \geq \vep_a, \; |\wh{\bs{p}} - p | \geq \vep_r p \} & = &
\Pr \{ |\wh{\bs{p}} - p | \geq \vep_r p \}\\
& = & \Pr \{ \wh{\bs{p}}  \geq p + \vep_r p \} + \Pr \{ \wh{\bs{p}}
\leq p - \vep_r p \}. \eee  Invoking Lemma \ref{lem6}, we have
\[
\Pr \{ \wh{\bs{p}}  \geq p + \vep_r p \} \leq \exp \li (n \; g \li (
\vep_a, \f{\vep_a}{\vep_r} \ri ) \ri ).
\]
On the other hand,
\[
 \Pr \{ \wh{\bs{p}}
\leq p - \vep_r p \} \leq \exp (n \; g( - \vep_r p, p)) \leq \exp
\li (n \; g \li ( - \vep_a, \f{\vep_a}{\vep_r}\ri ) \ri ) \leq \exp
\li (n \; g \li ( \vep_a, \f{\vep_a}{\vep_r}\ri ) \ri )
\]
where the first inequality follows from Lemma \ref{lem1}, the second inequality follows from Lemma \ref{lem4}, and the last inequality follows
from Lemma \ref{lem3}.  Hence,
\[
\Pr \{ |\wh{\bs{p}} - p | \geq \vep_a, \; |\wh{\bs{p}} - p | \geq
\vep_r p \}  \leq  2 \exp \li (n \; g \li ( \vep_a,
\f{\vep_a}{\vep_r} \ri ) \ri ) < \de. \] This proves (\ref{cov}) for
$\f{\vep_a}{\vep_r} < p < 1$. The proof of Theorem 1 is thus
completed.

\sect{Proof Theorem 2}

We need some preliminary results.  We shall introduce functions
\[
\mscr{M} (z, p) =  \ln \li ( \f{p}{z} \ri ) + \li ( \f{1}{z} - 1 \ri
) \ln \li (  \f{ 1 - p } { 1 - z } \ri ) \] and
\[
\mscr{H} (z, p) = z \; \mscr{M} (z, p)
\]
 for $0 < z < 1$ and $0 < p < 1$.

\beL Suppose $1 \leq r \leq M < N$.  Then,
\[
\Pr \li \{ \f{r}{\bs{n}} \leq (1 - \vep) p \ri \} \leq (1 - \vep)^{-
r} \exp \li ( - \f{ \vep r }{ 1 - \vep }  \ri ).
\]
\eeL

\bpf

Clearly, \bee \Pr \li \{ \f{r}{\bs{n}} \leq (1 - \vep) p \ri \} & = & \Pr \li \{ \bs{n} \geq \f{r}{(1 - \vep) p} \ri \}\\
& = & \Pr \{ \bs{n} \geq m \} \eee where
\[
m = \li  \lc \f{r}{(1 - \vep) p} \ri \rc.
\]
It can be seen that there exists a real number $\vep^* \in (0, 1)$
such that $\vep^* \geq \vep$ and
\[
\f{r}{(1 - \vep^*) p} = \li  \lc \f{r}{(1 - \vep) p} \ri \rc.
\]
Now let $K_m$ be the number of units having a certain attribute
among $m$ units drawn by a sampling without replacement from a
finite population of size $N$ with $M$ units having the attribute.
Then,  \bee \Pr \{ \bs{n} \geq m \} & = & \Pr \{  K_m \leq r \} \\
& = & \Pr \li \{ \f{K_m}{m} \leq \f{r}{m}   \ri \} \\
& = & \Pr \li \{ \f{K_m}{m} \leq (1 - \vep^*) p   \ri \}. \eee
Applying the well-known Hoeffding inequality \cite{Hoeffding} for
the case of finite population, we have \bee \Pr \li \{ \f{K_m}{m}
\leq (1 - \vep^*) p \ri \}
& \leq & \exp \li ( m  \mscr{H} (p - \vep^* p , p ) \ri )\\
& = & \exp \li (  \f{r}{(1 - \vep^*) p} \mscr{H} (p - \vep^* p , p ) \ri )\\
& = & \exp \li (   r   \mscr{M} (p - \vep^* p , p ) \ri )\\
& \leq &  \exp \li (   r   \mscr{M} (p - \vep p , p ) \ri ) \eee
where the last inequality follows from $\vep^* \geq \vep$ and the
monotone property of $\mscr{M} (p - \vep p , p )$ with respect to
$\vep$, which has been established as Lemma 5 in \cite{Chen3}.

From the proof of Lemma 6 of \cite{Chen3}, we know that $\mscr{M} (p
- \vep p , p )$ is monotonically decreasing with respect to $p \in
(0, 1)$. Hence,
\[
\Pr \li \{ \f{r}{\bs{n}} \leq (1 - \vep) p \ri \}  \leq \exp \li ( r
\mscr{M} (p - \vep p , p ) \ri ) \leq \lim_{p \to 0} \exp \li ( r
\mscr{M} (p - \vep p , p ) \ri ) = (1 - \vep)^{- r} \exp \li ( - \f{
\vep r }{ 1 - \vep }  \ri ).
\]
The proof of the lemma is thus completed.

 \epf

\beL Suppose $1 \leq r \leq M < N$ and $p + \vep p < 1$.  Then,
\[
\Pr \li \{ \f{r}{\bs{n}} \geq (1 + \vep) p \ri \} \leq (1 + \vep)^{-
r} \exp \li ( \f{ \vep r }{ 1 + \vep } \ri ).
\]
\eeL

\bpf

It is clear that \bee \Pr \li \{ \f{r}{\bs{n}} \geq (1 + \vep) p \ri \} & = & \Pr \li \{ \bs{n} \leq \f{r}{(1 + \vep) p} \ri \}\\
& = & \Pr \{ \bs{n} \leq m \} \eee where
\[
m =  \li  \lf \f{r}{(1 + \vep) p} \ri \rf.
\]
It can be seen that there exists a real number $\vep^* \in (0, 1)$
such that $\vep^* \geq \vep$ and
\[
\f{r}{(1 + \vep^*) p} = \li  \lf \f{r}{(1 + \vep) p} \ri \rf.
\]
Now let $K_m$ be the number of units having a certain attribute
among $m$ units drawn by a sampling without replacement from a
finite population of size $N$ with $M$ units having the attribute.
Then, \bee \Pr \{ \bs{n} \leq m \} & = & \Pr \{  K_m \geq r \} \\
& = & \Pr \li \{ \f{K_m}{m} \geq \f{r}{m}   \ri \} \\
& = & \Pr \li \{ \f{K_m}{m} \geq (1 + \vep^*) p   \ri \}. \eee
Applying the well-known Hoeffding inequality \cite{Hoeffding} for
the case of finite population, we have \bee \Pr \li \{ \f{K_m}{m}
\geq (1 + \vep^*) p \ri \} & \leq & \exp \li ( m  \mscr{H} ( \vep^* p , p ) \ri )\\
& = & \exp \li (  \f{r}{(1 + \vep^*) p} \mscr{H} ( \vep^* p , p ) \ri )\\
& = & \exp \li (   r   \mscr{M} (p + \vep^* p , p ) \ri )\\
& \leq &  \exp \li (   r   \mscr{M} (p + \vep p , p ) \ri ) \eee
where the last inequality follows from $\vep^* \geq \vep$ and the
monotone property of $\mscr{M} (p + \vep p , p )$ with respect to
$\vep$, which has been established as Lemma 5 in \cite{Chen3}.

From the proof of Lemma 6 of \cite{Chen3}, we know that $\mscr{M} (p
+ \vep p , p )$ is monotonically decreasing with respect to $p \in
\li ( 0, \f{1}{1 + \vep} \ri )$. Hence,
\[
\Pr \li \{ \f{r}{\bs{n}} \geq (1 + \vep) p \ri \}  \leq \exp \li ( r
\mscr{M} (p + \vep p , p ) \ri ) \leq \lim_{p \to 0} \exp \li ( r
\mscr{M} (p + \vep p , p ) \ri ) = (1 + \vep)^{- r} \exp \li (  \f{
\vep r }{ 1 + \vep }  \ri ).
\]
The proof of the lemma is thus completed.

\epf

\bsk

Now we are in a position to prove Theorem 2.  We shall consider the
following cases:

Case (i): $M < r$;

Case (ii): $M = N$;

Case (iii): $r = N$;

Case (iv): $1 \leq r \leq M < N$ and $p < \f{1}{1 + \vep}$;

Case (v): $1 \leq  r \leq M  < N$ and $p = \f{1}{1 + \vep}$;

Case (vi): $1 \leq r \leq M < N$ and $p > \f{1}{1 + \vep}$.

\bsk

In Case (i), we have $\bs{n} = N$ and $\bs{k} = M$.  Hence,
$\wt{\bs{p}} = p$ and $\Pr \li \{  \li |  \wt{\bs{p}} - p  \ri |
\geq \vep p \ri \} = 0 \leq \mscr{Q} (\vep, r)$.

\bsk

In Case (ii), we have $\wt{\bs{p}} = p$ and $\Pr \li \{  \li |
\wt{\bs{p}} - p  \ri | \geq \vep p \ri \} = 0 \leq \mscr{Q} (\vep,
r)$.

\bsk

In Case (iii), we have $\wt{\bs{p}} = p$ and $\Pr \li \{  \li |
\wt{\bs{p}} - p  \ri | \geq \vep p \ri \} = 0 \leq \mscr{Q} (\vep,
r)$.

\bsk

In Case (iv), we have $\bs{k} = r$ and, by Lemma 7 and Lemma 8, \bee
\Pr \li \{ \li | \wt{\bs{p}} - p  \ri | \geq \vep p \ri \} & = & \Pr
\li \{ \f{r}{\bs{n}} \leq (1 - \vep) p \ri \} + \Pr \li
\{ \f{r}{\bs{n}} \geq (1 + \vep) p \ri \} \\
& \leq &  (1 - \vep)^{- r} \exp \li ( - \f{ \vep r }{ 1 - \vep } \ri
) + (1 + \vep)^{- r} \exp \li ( \f{ \vep r }{ 1 + \vep } \ri )\\
& = & \mscr{Q} (\vep, r). \eee

\bsk

In Case (v), we have $\bs{k} = r$ and \bee \Pr \li \{  \li |
\wt{\bs{p}} - p  \ri | \geq \vep p \ri \} & = & \Pr \li \{
\f{r}{\bs{n}} \leq (1 - \vep) p \ri \} + \Pr \li
\{ \f{r}{\bs{n}} \geq (1 + \vep) p \ri \} \\
& = & \Pr \li \{ \f{r}{\bs{n}} \leq (1 - \vep) p \ri \} + \Pr \li \{
\bs{k} = \bs{n} = r \ri \}. \eee Notice that \bee \Pr \li \{ \bs{k}
= \bs{n} = r \ri \} & = & \f{ \bi{M}{ r} } { \bi{N}{r} } < \li (
\f{M}{N} \ri )^r =  p^r = \li ( \f{1}{1 + \vep} \ri )^r < (1 +
\vep)^{- r} \exp \li ( \f{ \vep r }{ 1 + \vep } \ri ) \eee as a
result of $M \leq N$.  Therefore, by Lemma 7,
\[
\Pr \li \{  \li | \wt{\bs{p}} - p  \ri | \geq \vep p \ri \} \leq (1
- \vep)^{- r} \exp \li ( - \f{ \vep r }{ 1 - \vep } \ri ) + (1 +
\vep)^{- r} \exp \li ( \f{ \vep r }{ 1 + \vep } \ri ) = \mscr{Q}
(\vep, r).
\]

\bsk

In Case (vi), we have $\bs{k} = r, \; \Pr \li \{ \f{r}{\bs{n}} \geq
(1 + \vep) p \ri \} = 0$ and, by Lemma 7, \bee \Pr \li \{  \li |
\wt{\bs{p}} - p \ri | \geq \vep p \ri \} & = & \Pr \li \{
\f{r}{\bs{n}} \leq (1 - \vep) p \ri \} + \Pr \li
\{ \f{r}{\bs{n}} \geq (1 + \vep) p \ri \} \\
& = & \Pr \li \{ \f{r}{\bs{n}} \leq (1 - \vep) p \ri \}\\
 & \leq & (1
- \vep)^{- r} \exp \li ( - \f{ \vep r }{ 1 - \vep } \ri )  <
\mscr{Q} (\vep, r). \eee

\bsk

So, we have shown $\Pr \li \{  \li | \wt{\bs{p}} - p \ri | \geq \vep
p \ri \} \leq \mscr{Q} (\vep, r)$.  The other statements of Theorem
2 have been established in \cite{Chen3}.

This concludes the proof of Theorem 2.

\end{document}